\documentclass[preprint,11pt,compress]{elsarticle}

\usepackage{amssymb}
\usepackage{amsmath}
\usepackage{bbm}
\usepackage[a4paper, body={15cm,23cm}]{geometry}

\begin{document}
\journal{arXiv}
\parindent 15pt
\parskip 5pt
 \newcommand{\RE}{\forall}
 \newcommand{\eps}{\varepsilon}
 \newcommand{\lam}{\lambda}
 \newcommand{\To}{\rightarrow}
 \newcommand{\as}{{\rm d}\mathbbm{P}\times {\rm d}t-a.e.}
 \newcommand{\ps}{{\rm d}\mathbbm{P}-a.s.}
 \newcommand{\jf}{\int_t^T}
 \newcommand{\tim}{\times}
 \newcommand{\T}{[0,T]}
 \newcommand{\Lp}{L^p(\R^k)}

 \newcommand{\D}{\mathcal{D}} \newcommand{\F}{\mathbbm{F}}
 \newcommand{\E}{\mathbbm{E}}
 \newcommand{\N}{\mathbbm{N}}
 \newcommand{\s}{\mathcal{S}}
 \newcommand{\R}{\mathbbm{R}}

 \newcommand {\Dis}{\displaystyle}


\begin{frontmatter}

\title {$L^p$ solutions of multidimensional BSDEs with
 weak monotonicity and general growth generators\vspace{-0.2cm}\tnoteref{fund}}
\tnotetext[fund]{Supported by the National Natural Science Foundation of China (No. 11371362), the China Postdoctoral Science Foundation (No. 2013M530173), the Qing Lan Project and the Fundamental Research Funds for the Central Universities (No. 2013RC20)\vspace{0.2cm}.}

\author{ShengJun FAN}
\ead{f$\_$s$\_$j@126.com}

\address{College of Sciences, China University of Mining and Technology, Xuzhou 221116, PR China\vspace{-1.0cm}}

\begin{abstract}
In this paper, we first establish the existence and uniqueness of $L^p\ (p>1)$ solutions for multidimensional backward stochastic differential equations (BSDEs) under a weak monotonicity condition together with a general growth condition in $y$ for the generator $g$. Then, we overview several conditions related closely to the weak monotonicity condition and compare them in an effective way. Finally, we put forward and prove a stability theorem and a comparison theorem of $L^p\ (p>1)$ solutions for this kind of BSDEs.
\end{abstract}

\begin{keyword}
Backward stochastic differential equation\sep $L^p$ solution \sep Weak monotonicity condition\sep Comparison theorem \sep Stability theorem  \vspace{0.2cm}

\MSC[2010] 60H10
\end{keyword}
\end{frontmatter}\vspace{-0.4cm}


\section{Introduction}

Throughout this paper, let us fix a real number $T>0$, and two positive integers $k$ and $d$. Let $\R^+:=[0,+\infty)$ and let $(\Omega,\F,\mathbbm{P})$ be a probability space carrying a standard $d$-dimensional Brownian motion $(B_t)_{t\geq
0}$, $(\F_t)_{t\geq 0}$ be the natural $\sigma$-algebra generated by $(B_t)_{t\geq 0}$ and $\F=\F_T$. For each
subset $A\subset \Omega\times \T$, let $\mathbbm{1}_{A}=1$ in case of $(t,\omega)\in A$, otherwise, let $\mathbbm{1}_{A}=0$. The Euclidean norm of a vector $y\in\R^k$ will be defined by $|y|$, and for an $k\times d$ matrix $z$, we define $|z|=\sqrt{{\rm Tr}zz^*}$,
where $z^*$ is the transpose of $z$. Let $\langle x,y\rangle$ represent the inner product of $x,y\in\R^k$. For each $p>1$, we denote by $\Lp$ the set of all $\R^k$-valued and $\F_T$-measurable random vectors $\xi$ such that $\E[|\xi|^p]<+\infty$, and by ${\s}^p(0,T;\R^k)$ the set of $\R^k$-valued, adapted and continuous processes $(Y_t)_{t\in\T}$ such that
$$\|Y\|_{{\s}^p}:=\left( \E[\sup_{t\in\T} |Y_t|^p] \right)
^{1/p}<+\infty.$$
Moreover, let ${\rm M}^p(0,T;\R^{k\times d})$ denote the set of $(\F_t)$-progressively
measurable ${\R}^{k\times d}$-valued processes $(Z_t)_{ t\in\T}$ such that
$$\|Z\|_{{\rm M}^p}:=\left\{ \E\left[\left(\int_0^T |Z_t|^2\
{\rm d}t\right)^{p/2}\right] \right\}^{1/p}<+\infty.
$$
Obviously, both ${\s}^p$ and ${\rm M}^p$ are Banach spaces for each $p>1$.

In this paper, we are concerned with the following multidimensional backward stochastic differential equation (BSDE for short in the remaining):
\begin{equation}
    y_t=\xi+\int_t^Tg(s,y_s,z_s){\rm d}s-\int_t^Tz_s {\rm d}B_s,\ \
    t\in\T,
\end{equation}
where $\xi\in\Lp$ is called the terminal condition, $T$ is called the time horizon, the random function $$g(\omega,t,y,z):\Omega\tim \T\tim {\R}^{k }\tim
{\R}^{k\times d}\longmapsto {\R}^k$$
is $(\F_t)$-progressively measurable for each $(y,z)$, called the generator of BSDE (1). This BSDE is usually denoted by the BSDE ($\xi,T,g$).

For convenience of the following discussion, we introduce the following definitions concerning solutions of BSDE (1).

{\bf Definition 1}\ \ A solution to BSDE (1) is a pair of $(\F_t)$-progressively measurable processes $(y_t,z_t)_{t\in\T}$ with values in ${\R}^k\times {\R}^{k\times d}$ such that $\ps$, $t\mapsto y_t$ is continuous, $t\mapsto z_t$ belongs to ${\rm L}^2(0,T)$, $t\mapsto g(t,y_t,z_t)$ belongs to ${\rm L}^1(0,T)$, and $\ps$, (1) holds true for each $t\in\T$.\vspace{0.2cm}

{\bf Definition 2}\ \ Assume that $(y_t,z_t)$ is a solution to BSDE (1). If $(y_t,z_t)\in {\s}^p(0,T;\R^{k})\times {\rm M}^p(0,T;\R^{k\times d})$ for some $p>1$, then it will be called an $L^p$ solution of BSDE (1).

Nonlinear BSDEs were firstly introduced in 1990 by \citet{Par90}, who established the existence and uniqueness for $L^2$ solutions of BSDEs under the Lipschitz assumption of the generator $g$. Since then, BSDEs have been studied with great interest, and they have become a powerful tool in many fields above all financial mathematics, stochastic games and optimal control, non-linear PDEs and homogenization. See \cite{BC08,Buc99,Buc00,Delb11,Delb10,El97,Kob00,Mor09,
Par99,Peng91,Peng97,Tang94,Tang98}
and the references therein for applications of BSDEs to PDEs, optimal control, homogenization as well as in mathematical finances.

From the beginning, many authors attempted to improve the result of \cite{Par90} by weakening the Lipschitz hypothesis on $g$, see \cite{Bah02,Bah10,BC08,Bri06,Bri08,Bri07,Cons01,Delb11,
Delb10,El97,Fan10,FJ11,Fan12,FJ12a,Fan13,Fan14,Fan10b,Ham03,Jia08,
Jia10,Kob00,Lep97,Lep98,Mao95,Mor09,Par99,Peng91,Situ97,Xu14,
Yin08,Yin03}, or the $L^2$ integrability assumptions on $\xi$, see \cite{Bri03,Bri06,Bri08,Chen10,FJ12b,Fan14,
Peng97,Yao10}, or relaxing the finite terminal time $T$ to a stopping time or infinity, see \cite{Chen00,Fan11,Kob00,Mor09,Par99,Yin08}. From these results we can see that the case of one-dimensional BSDEs is easier to handle due to the presence of the comparison theorem of solutions (see \cite{Bri06,Bri08,Bri07,Chen10,Chen00,Fan10,FJ11,
Fan12,FJ12a,FJ12b,Fan11,Jia08,Jia10,Kob00,Lep97,Lep98,
Peng97}).

One of the main purposes of the present paper is to establish an existence and uniqueness of $L^p\ (p>1)$ solutions for multidimensional BSDEs under weaker conditions on the generators. Here, we would like to mention the following several results on multidimensional BSDEs, which is related closely to our result. First of all, \citet{Mao95} obtained an existence and uniqueness result of an $L^2$ solution for (1) where $g$ satisfies a particular non-Lipschitz condition in $y$ called usually the Mao condition in the literature, and \citet{Fan14} investigated the existence and uniqueness of an $L^p\
(p>1)$ solution for (1) where $g$ satisfies a new kind of non-Lipschitz condition in $y$. Second, \citet{Peng91} first introduced a kind of monotonicity condition in $y$ for $g$, and under this monotonicity condition as well as a general growth condition in $y$ for $g$, \citet{Par99} established an existence and uniqueness result of an $L^2$ solution for (1). Using the same monotonicity condition and a more general growth condition in $y$ for $g$, \citet{Bri03} investigated the existence and uniqueness of an $L^p\ (p\geq 1)$ solution for (1). Furthermore, \citet{Situ97} put forward a kind of weak monotonicity condition in $y$ for $g$ and considered the existence and uniqueness of $L^p\ (p\geq 1)$ solutions for BSDEs with jumps, but the generator $g$ is forced to also satisfy a linear growth condition in $y$. Recently, \citet{Fan13} and \citet{Xu14} established the existence and uniqueness of an $L^2$ solution for (1) under the weak monotonicity condition and the more general growth condition in $y$ for the generator $g$, which really and truly unifies the Mao condition in $y$ and the monotonicity condition with the general growth condition in $y$.

In this paper, we first establish the existence and uniqueness of $L^p\ (p>1)$ solutions for multidimensional BSDEs under the weak monotonicity condition together with the more general growth condition in $y$ for the generator $g$ (see Theorem 1 in Section 2 and its proof in Section 4), which extends some existing results including Theorem 4.2 in \citet{Bri03} and Theorem 1 in \citet{Fan14}. Then, we overview several conditions related closely to the weak monotonicity condition and compare them in an effective way (see Proposition 1 in Section 2 and its proof in Appendix). Finally, we put forward and prove a stability theorem and a comparison theorem of $L^p\ (p>1)$ solutions for this kind of BSDEs (see Theorem 2 in Section 4 and Theorem 3 in Section 5).

This paper is organized as follows. In Section 2 we state the assumptions and the existence and uniqueness result for $L^p\ (p>1)$ solutions of multidimensional BSDEs and introduce several propositions, corollaries, remarks and examples to show that it generalizes some existing results. In Section 3, we establish two nonstandard a priori estimates for $L^p\ (p>1)$ solutions of multidimensional BSDEs, based on which we prove a stability Theorem and the existence and uniqueness result in Section 4. Then, we put forward and prove a new comparison theorem for $L^p\ (p>1)$ solutions of one dimensional BSDEs in Section 5. Finally, the proof of the relations between the assumptions related closely to the weak monotonicity condition is provided in Appendix.

\section{An existence and uniqueness result}

In this section, we will state the existence and
uniqueness result for $L^p\ (p>1)$ solutions of multidimensional BSDEs and introduce several propositions, corollaries, remarks and examples to show that it generalizes some existing results including Theorem 4.2 in \citet{Bri03} and Theorem 1 in \citet{Fan14}. Let us start with  introducing the following assumptions:\vspace{0.2cm}

{\bf (H1)$_p$} $g$ satisfies the $p$-order weak monotonicity condition in $y$, i.e., there exists a nondecreasing and concave function $\rho(\cdot):\R^+\mapsto \R^+$ with $\rho(0)=0$, $\rho(u)>0$ for $u>0$ and $\int_{0^+} {{\rm d}u\over \rho(u)}=+\infty$ such that $\as$, $\RE\ y_1,y_2\in \R^k,z\in\R^{k\times d},$\vspace{0.25cm}\\
\hspace*{1.5cm}$|y_1-y_2|^{p-1}\langle {y_1-y_2\over |y_1-y_2|}\mathbbm{1}_{|y_1-y_2|\neq 0},g(\omega,t,y_1,z)-g(\omega,t,y_2,z)\rangle\leq \rho(|y_1-y_2|^p);$\vspace{0.3cm}

{\bf (H2)} $\as$, $\RE\ z\in {\R^{k\times d}},\ \ \ y\longmapsto g(\omega,t,y,z)$ is continuous;\vspace{0.25cm}

{\bf (H3)} $g$ has a general growth with respect to $y$, i.e,
$$\RE\ \ \alpha>0,\ \phi_\alpha(t):=\sup\limits_{|y|\leq \alpha}
|g(\omega,t,y,0)-g(\omega,t,0,0)|\in L^1(\T\tim \Omega);$$

{\bf (H4)}\ $g$ is Lipschitz continuous in $z$, uniformly with respect to $(\omega,t,y)$, i.e., there exists a
constant $\bar \lam\geq 0$ such that $\as$, $\RE\ y\in \R^k,z_1,z_2\in\R^{k\times d}$,
$$|g(\omega,t,y,z_1)-g(\omega,t,y,z_2)|\leq \bar \lam |z_1-z_2|;$$

{\bf (H5)$_p$}\ $\Dis\E\left[|\xi|^p+\left(\int_0^T |g(\omega,t,0,0)|\ {\rm d}t\right)^p\right]<+\infty$.\vspace{0.25cm}

The following Theorem 1 is one of the main results of this paper. It's proof will be given in Section 4.

{\bf Theorem 1}\ \ Assume that $p>1$, and assumptions (H1)$_{p\wedge 2}$, (H2)-(H4) and (H5)$_p$ hold. Then, the BSDE $(\xi,T,g)$ has a unique $L^p$ solution.\vspace{0.2cm}

It should be mentioned that Theorem 1 has been proved in \citet{Xu14} for the case of $p=2$. In addition, by Theorem 1 the following corollary is immediate.

{\bf Corollary 1}\ \ Assume that the generator $g$ satisfies assumptions (H1)$_{2}$ and (H2)-(H4). Then, if (H5)$_p$ holds for some $p>2$, then the BSDE $(\xi,T,g)$ has a unique $L^p$ solution.\vspace{0.2cm}

In the sequel, let us further introduce the following assumptions on $g$:

{\bf (H1a)$_p$} $g$ satisfies the $p$-order one-sided Mao condition in $y$, i.e., there exists a nondecreasing and concave function $\rho(\cdot):\R^+\mapsto \R^+$ with $\rho(0)=0$, $\rho(u)>0$ for $u>0$ and $\int_{0^+} {{\rm d}u\over \rho(u)}=+\infty$ such that $\as$, $\RE\ y_1,y_2\in \R^k,z\in\R^{k\times d},$\vspace{0.25cm}\\
\hspace*{1.5cm}$\langle {y_1-y_2\over |y_1-y_2|}\mathbbm{1}_{|y_1-y_2|\neq 0},
g(\omega,t,y_1,z)-g(\omega,t,y_2,z)\rangle\leq \rho^{{1\over p}}(|y_1-y_2|^p);$\vspace{0.25cm}

{\bf (H1b)$_p$} $g$ satisfies the $p$-order one-sided Costantin condition in $y$, i.e., there exists a nondecreasing and concave function $\rho(\cdot):\R^+\mapsto \R^+$ with $\rho(0)=0$, $\rho(u)>0$ for $u>0$ and $\int_{0^+} {u^{p-1}\over \rho^p(u)}{\rm d}u=+\infty$ such that $\as$, $\RE\ y_1,y_2\in \R^k,z\in\R^{k\times d},$\vspace{0.25cm}\\
\hspace*{1.5cm}$\langle {y_1-y_2\over |y_1-y_2|}\mathbbm{1}_{|y_1-y_2|\neq 0},
g(\omega,t,y_1,z)-g(\omega,t,y_2,z)\rangle\leq \rho(|y_1-y_2|);$\vspace{0.25cm}

{\bf (H1*)} $g$ satisfies the one-sided Osgood condition in $y$, i.e., there exists a nondecreasing and concave function $\rho(\cdot):\R^+\mapsto \R^+$ with $\rho(0)=0$, $\rho(u)>0$ for $u>0$ and $\int_{0^+} {{\rm d}u\over \rho(u)}=+\infty$ such that $\as$, $\RE\ y_1,y_2\in \R^k,z\in\R^{k\times d},$\vspace{0.25cm}\\
\hspace*{1.5cm}$\langle {y_1-y_2\over |y_1-y_2|}\mathbbm{1}_{|y_1-y_2|\neq 0},
g(\omega,t,y_1,z)-g(\omega,t,y_2,z)\rangle\leq \rho(|y_1-y_2|).$\vspace{0.55cm}

{\bf Remark 1}\  It is easy to see that the following statements are true.
\begin{itemize}
\item When $\rho(x)=\mu x$ for some constant $\mu>0$, (H1)$_p$, (H1a)$_p$, (H1b)$_p$ and (H1*) are all the known monotonicity condition for each $p\geq 1$;
\item In case of $p=1$, (H1)$_p$, (H1a)$_p$ and (H1b)$_p$ are all the same as (H1*);
\item In case of $p=2$, (H1)$_p$, (H1a)$_p$ and (H1b)$_p$ are respectively the so-called weak monotonicity condition, one-sided Mao condition and one-sided Constantin condition put forward in \citet{Fan13}.\vspace{0.2cm}
\end{itemize}

With respect to the previous assumptions, we have the following important observation. It's proof will be provided in Appendix.

{\bf Proposition 1} For each $1\leq p\leq q<+\infty$, we have\vspace{0.1cm}
$$
\begin{array}{cccccc}
(i)& {\rm (H1^*)} &\Longrightarrow &{\rm (H1)}_p &\Longrightarrow &{\rm (H1)}_q;\\
(ii)& {\rm (H1b)}_q& \Longrightarrow& {\rm (H1b)}_p& \Longrightarrow& {\rm (H1^*)};\\
(iii)& {\rm (H1a)}_p& \Longleftrightarrow& {\rm (H1b)}_p.&
\end{array}$$
In addition, we can show that for each $p\geq 1$, the concavity condition of $\rho(\cdot)$ in assumptions (H1a)$_p$ and (H1b)$_p$ can be replaced with the continuity condition.\vspace{0.2cm}

According to Theorem 1 and Proposition 1, the following corollaries follow immediately.

{\bf Corollary 2}\ \ Assume that the generator $g$ satisfies assumptions (H1*) and (H2)-(H4). Then, if (H5)$_p$ holds for some $p>1$, then the BSDE $(\xi,T,g)$ has a unique $L^p$ solution.\vspace{0.2cm}

{\bf Corollary 3}\ \ Assume that $p>1$, and assumptions (H1a)$_{p}$ (or (H1b)$_{p}$), (H2)-(H4) and (H5)$_p$ hold. Then, the BSDE $(\xi,T,g)$ has a unique $L^p$ solution.\vspace{0.3cm}

The following four assumptions (H1')$_p$, (H1a')$_p$, (H1b')$_p$ and (H1'*) are respectively the stronger and two-sided versions of assumptions (H1)$_p$, (H1a)$_p$, (H1b)$_p$ and (H1*): \vspace{0.2cm}

{\bf (H1')$_p$} There exists a nondecreasing and concave function $\rho(\cdot):\R^+\mapsto \R^+$ with $\rho(0)=0$, $\rho(u)>0$ for $u>0$ and $\int_{0^+} {{\rm d}u\over \rho(u)}=+\infty$ such that $\as$, $\RE\ y_1,y_2\in \R^k,z\in\R^{k\times d},$\vspace{0.25cm}\\
\hspace*{1.5cm}$|y_1-y_2|^{p-1}|g(\omega,t,y_1,z)-
g(\omega,t,y_2,z)|\leq \rho(|y_1-y_2|^p);$\vspace{0.25cm}

{\bf (H1a')$_p$} $g$ satisfies the $p$-order Mao condition in $y$, i.e., there exists a nondecreasing and concave function $\rho(\cdot):\R^+\mapsto \R^+$ with $\rho(0)=0$, $\rho(u)>0$ for $u>0$ and $\int_{0^+} {{\rm d}u\over \rho(u)}=+\infty$ such that $\as$, $\RE\ y_1,y_2\in \R^k,z\in\R^{k\times d},$\vspace{0.25cm}\\
\hspace*{1.5cm}$|g(\omega,t,y_1,z)-g(\omega,t,y_2,z)|\leq \rho^{{1\over p}}(|y_1-y_2|^p);$\vspace{0.25cm}

{\bf (H1b')$_p$} $g$ satisfies the $p$-order Constantin condition in $y$, i.e., there exists a nondecreasing and concave function $\rho(\cdot):\R^+\mapsto \R^+$ with $\rho(0)=0$, $\rho(u)>0$ for $u>0$ and $\int_{0^+} {u^{p-1}\over \rho^p(u)}{\rm d}u=+\infty$ such that $\as$, $\RE\ y_1,y_2\in \R^k,z\in\R^{k\times d},$\vspace{0.25cm}\\
\hspace*{1.5cm}$|g(\omega,t,y_1,z)-g(\omega,t,y_2,z)|\leq \rho(|y_1-y_2|);$\vspace{0.25cm}

{\bf (H1'*)} $g$ satisfies the Osgood condition in $y$, i.e., there exists a nondecreasing and concave function $\rho(\cdot):\R^+\mapsto \R^+$ with $\rho(0)=0$, $\rho(u)>0$ for $u>0$ and $\int_{0^+} {{\rm d}u\over \rho(u)}=+\infty$ such that $\as$, $\RE\ y_1,y_2\in \R^k,z\in\R^{k\times d},$\vspace{0.25cm}\\
\hspace*{1.5cm}$|g(\omega,t,y_1,z)-g(\omega,t,y_2,z)|\leq \rho(|y_1-y_2|).$\vspace{0.55cm}

{\bf Remark 2} It is easy to see that the following statements are true.
\begin{itemize}
\item When $\rho(x)=\mu x$ for some constant $\mu>0$, (H1')$_p$, (H1a')$_p$, (H1b')$_p$ and (H1'*) are all the known Lipschitz condition for each $p\geq 1$;
\item In case of $p=1$, (H1')$_p$, (H1a')$_p$ and (H1b')$_p$ are the same as (H1'*);
\item In case of $p=2$, (H1a')$_p$ and (H1b')$_p$ are respectively the known Mao condition and Constantin condition;
\item For each $p\geq 1$, we have (H1')$_p$$\Longrightarrow$(H1)$_p$, (H1a')$_p$$\Longrightarrow$(H1a)$_p$$+$(H2)$+$(H3), and\\ (H1b')$_p$$\Longrightarrow$(H1b)$_p$$+$(H2)$+$(H3);
\item Proposition 1 holds also true for assumptions (H1')$_p$, (H1a')$_p$, (H1b')$_p$ and (H1*).
\end{itemize}

{\bf Remark 3}\ It follows from Remarks 1-2 and Proposition 1 that Theorem 1 and some of its corollaries all improve some existing results for $L^p$ solutions of multidimensional BSDEs including Theorem 4.2 in \citet{Bri03} and Theorem 1 in \citet{Fan14}.\vspace{0.2cm}

Now, we give two examples of BSDEs which satisfy the assumptions in Corollary 2. In our knowledge, they are not covered by the previous works.

{\bf Example 1} Let $k=1$, $\bar p\geq 1$ and
$$g(\omega,t,y,z)=h(|y|)-e^{|B_t(\omega)|\cdot y}+(e^{-y}\wedge 1)\cdot |z|+{1\over \sqrt[3]{t}}\mathbbm{1}_{t>0},$$
where
$$h(x)=\left\{
\begin{array}{lll}
-x|\ln x|^{1/\bar p}& ,&0<x\leq \delta;\\
h'(\delta-)(x-\delta)+h(\delta)& ,&x> \delta;\\
0& ,&{\rm other\ cases}.
\end{array}\right.\vspace{0.1cm}$$
with $\delta>0$ small enough.

It is easy to see that $g$ satisfies (H2)-(H4) with $\bar\lam =1$. Furthermore, we can prove that $g$ satisfies (H1b)$_{\bar p}$ by verifying that $e^{-\beta y}$ with $\beta\geq 0$ is decreasing in $y$, $h(\cdot)$ is concave and sub-additive on $\R^+$ and then the following inequality holds: $\as$,
$$\RE\ y_1,y_2,z,\ \
\langle {y_1-y_2\over |y_1-y_2|}\mathbbm{1}_{|y_1-y_2|\neq 0},
g(\omega,t,y_1,z)-g(\omega,t,y_2,z)\rangle\leq
h(|y_1-y_2|)$$
with
$$\int_{0^+} {u^{\bar p-1}\over h^{\bar p}(u)}{\rm d}u=+\infty.\vspace{0.1cm}$$
Thus, (H1*) holds for $g$ by Proposition 1, and then from Corollary 2 we know that if $\xi\in\Lp$ for some $p>1$, then BSDE $(\xi,T,g)$ has a unique $L^p$ solution.\vspace{0.2cm}

{\bf Example 2}\ Let $y=(y_1,\cdots,y_k)$ and
$g(t,y,z)=(g_1(t,y,z),\cdots,g_k(t,y,z))$, where for each $i=1,\cdots,k$,
$$g_i(\omega,t,y,z):=e^{-y_i}+h(|y|)+\sin |z|+|B_t(\omega)|,$$
and $h(x)$ is defined in Example 1.

It is not hard to verify that this generator $g$ satisfies (H1b)$_{\bar p}$, (H1*) and (H2)-(H4) with $\bar\lam =1$. It then follows from Corollary 2 that if $\xi\in\Lp$ for some $p>1$, then BSDE $(\xi,T,g)$ has a unique $L^p$ solution.

Finally, we would like to mention that the function $h(x)$ defined in Example 1 satisfies that
$$\RE\ q>\bar p,\ \ \int_{0^+}{u^{q-1}\over h^q(u)}{\rm d}u<+\infty.$$
And, we can also prove that neither of the two generators $g$ defined in Examples 1-2 satisfies (H1a)$_q$ or (H1b)$_q$ for each $q>\bar p$, which means that the inverse version of $(ii)$ of Proposition 1 does not hold.

\section{Two nonstandard a priori estimates}

In this section, we will establish two nonstandard a priori estimates concerning $L^p$ solutions of multidimensional BSDE (1), which will play an important role in the proof of our main results. The following assumption on the generator $g$ will be used:\vspace{0.2cm}

{\bf (A1)}\ \ $\as,\ \RE\ (y,z)\in \R^k\times\R^{k\times d}$,\vspace{0.25cm}\\
\hspace*{4cm}$\langle y,
g(\omega,t,y,z)\rangle\leq \mu |y|^2+\lam |y||z|+|y|f_t+\varphi_t,\ \ $\vspace{0.3cm}\\
where $\mu$ and $\lam$ are two non-negative constants, $f_t$ and $\varphi_t$ are two non-negative and $(\F_t)$-progressively measurable processes with
$$
\E\left[\left(\int_0^T f_t\ {\rm d}t\right)^p\right]<+\infty\ \ \ {\rm and}\ \ \ \E\left[\left(\int_0^T \varphi_t\ {\rm d}t\right)^{p/2}\right]<+\infty.
$$

{\bf Proposition 2} Assume that $p>0$ and (A1) holds. Let $(y_t,z_t)_{t\in\T}$ be a solution of BSDE (1) such that $y_t$ belongs to ${\s}^p(0,T;\R^{k})$. Then $z_t$ belongs to ${\rm M}^p(0,T;\R^{k\times d})$, and for each $0\leq u\leq t\leq T$, we have
$$\begin{array}{lll}
\Dis \E\left[\left.\left(\int_t^T |z_s|^2\ {\rm
d}s\right)^{p/2}\right|\F_u\right]&\leq & \Dis
C_{\mu,\lam,p,T}\E\left[\left.\sup\limits_{s\in
[t,T]}|y_s|^p\right|\F_u\right]+C_p\E\left[\left.\left(\int_t^T f_s\ {\rm d}s\right)^p\right|\F_u\right]\\
&&\Dis +C_p\E\left[\left.\left(\int_t^T \varphi_s\ {\rm
d}s\right)^{p/2}\right|\F_u\right],
\end{array}$$
where $C_{\mu,\lam,p,T}$ is a nonnegative constant depending on $(\mu,\lam,p,T)$, and $C_p$ is a nonnegative constant depending only on $p$.\vspace{0.2cm}

{\bf Remark 4}\ \ Note that the constant $C_p$ does not depend on $\mu$ and $\lam$. This fact will play an important role later.\vspace{0.2cm}

{\bf Proof of Proposition 2.}\ For each integer $n\geq 1$, let us introduce the stopping time
$$\tau_n=\inf \left\{t\in\T: \int_0^t |z_s|^2\
{\rm d}s\geq n\right\}\wedge T.\vspace{0.2cm}$$
Applying
It\^{o}'s formula to $|y_t|^2$ leads the equation
$$
|y_{t\wedge \tau_n}|^2+\int_{t\wedge \tau_n}^{\tau_n} |z_s|^2\ {\rm
d}s=|y_{\tau_n}|^2+2\int_{t\wedge \tau_n}^{\tau_n} \langle
y_s,g(s,y_s,z_s)\rangle \ {\rm d}s-2\int_{t\wedge
\tau_n}^{\tau_n}\langle y_s,z_s{\rm d}B_s\rangle,\ t\in\T. $$
It follows from (A1) that for each $s\in [{t\wedge
\tau_n},{\tau_n}]$,
$$
2\langle y_s,g(s,y_s,z_s)\rangle \leq 2(\mu+\lambda^2) |y_s|^2+{|z_s|^2\over 2}+2|y_s|f_s+2\varphi_s.
$$
Thus, we have
$$
\begin{array}{lll}
\Dis {1\over 2}\int_{t\wedge \tau_n}^{\tau_n} |z_s|^2\ {\rm d}s &\leq & \Dis 2[(\mu+\lambda^2)T +1]\sup\limits_{s\in [{t\wedge
\tau_n},T]}|y_s|^2+\left(\int_{t\wedge \tau_n}^T f_s\ {\rm d}s\right)^2\\
&& \Dis \ +2\int_{t\wedge \tau_n}^T \varphi_s\ {\rm d}s+2\left|\ \int_{t\wedge \tau_n}^{\tau_n}\langle
y_s,z_s{\rm d}B_s\rangle\right|,
\end{array}
$$
and the inequality $(a+b)^{p/2}\leq 2^p(a^{p/2}+b^{p/2})$ yields the existence of a constant $c_p>0$ depending only on $p$ such that
\begin{equation}
\begin{array}{lll}
\Dis \left(\int_{t\wedge \tau_n}^{\tau_n} |z_s|^2\ {\rm d}s \right)^{p/2}&\leq & \Dis c_p [(\mu+\lambda^2)T +1]^{p/2}\sup\limits_{s\in [{t\wedge
\tau_n},T]}|y_s|^p+c_p\left(\int_{t\wedge \tau_n}^T f_s\ {\rm d}s\right)^p\\
&&\Dis +c_p\left(\int_{t\wedge \tau_n}^T \varphi_s\ {\rm d}s\right)^{p/2}+c_p\left|\ \int_{t\wedge \tau_n}^{\tau_n}\langle y_s,z_s{\rm d}B_s\rangle\right|^{p/2}.
\end{array}
\end{equation}
Furthermore, the Burkholder-Davis-Gundy (BDG) inequality yields that there exists a constant $d_p>0$ depending only on $p$ such that for each $0\leq u\leq t\leq T$,
$$
\begin{array}{ll}
&\Dis c_p\E\left[\left.\left|\ \int_{t\wedge \tau_n}^{\tau_n}\langle y_s,z_s{\rm d}B_s\rangle\right|^{p/2}\right|\F_u\right]\\
\leq & \Dis
d_p\E\left[\left.\left(\int_{t\wedge \tau_n}^{\tau_n}
|y_s|^2|z_s|^2\ {\rm d}s\right)^{p/4}
\right|\F_u\right]\\
\leq & \Dis {d_p^2\over 2}\E\left[\left.\sup\limits_{s\in [{t\wedge\tau_n},T]}|y_s|^p\right|\F_u\right]+{1\over 2}\E\left[\left. \left(\int_{t\wedge \tau_n}^{\tau_n}|z_s|^2\ {\rm d}s\right)^{p/2}\right|\F_u\right].
\end{array}
$$
Finally, taking the conditional mathematical expectation with respect to $\F_u$ in both sides of (2) and using the above inequality together with Fatou's lemma and Lebesgue's dominated convergence theorem yields the desired result. The proof is completed.\vspace{0.2cm}\hfill $\square$

Let us further introduce the following assumption on the generator $g$:\vspace{0.2cm}

{\bf (A2)}\ \ $\as,\ \RE\ (y,z)\in \R^k\times\R^{k\times d}$,\vspace{0.25cm}\\
\hspace*{4cm}$|y|^{p-1}\langle {y\over |y|}\mathbbm{1}_{|y|\neq 0},
g(\omega,t,y,z)\rangle\leq \psi(|y|^p)+\lam |y|^{p-1}|z|+|y|^{p-1}f_t,$\vspace{0.3cm}\\
where $\lam$ is a non-negative constant, $f_t$ is a non-negative and $(\F_t)$-progressively measurable process with
$$\E\left[\left(\int_0^T f_t\ {\rm d}t\right)^p\right]<+\infty,$$
and $\psi(\cdot):\R^+\mapsto \R^+$ is a nondecreasing and concave function with $\psi(0)=0$.\vspace{0.1cm}

{\bf Proposition 3} Assume that $p>1$ and (A2) holds. Let $(y_t,z_t)_{t\in\T}$ be an $L^p$ solution of BSDE (1). Then, there exists a nonnegative constant $C_{\lam,p,T}$ depending only on $\lambda$, $p$ and $T$ such that for each $0\leq u\leq t\leq T$,
$$
\begin{array}{lll}
&&\Dis \E\left[\left.\sup\limits_{s\in
[t,T]}|y_s|^{p}\right|\F_u\right]\\
&\leq & \Dis
C_{\lam,p,T}\left\{\E[\left.|\xi|^p\right|\F_u]+\int_t^T
\psi(\E[\left.|y_s|^p\right|\F_u]) \ {\rm
d}s+\E\left[\left.\left(\int_t^T f_s\ {\rm
d}s\right)^{p}\right|\F_u\right]\right\}.\vspace{0.1cm}
\end{array}
$$

{\bf Proof.}\ It follows from  Corollary 2.3 in \citet{Bri03} that, with $c(p)=p[(p-1)\wedge 1]/2$,
$$\begin{array}{lll}
&&\Dis |y_t|^p+c(p)\int_t^T |y_s|^{p-2}\mathbbm{1}_{|y_s|\neq 0}|z_s|^2\ {\rm d}s\\ &\leq & \Dis|\xi|^p +p\int_t^T |y_s|^{p-2}\mathbbm{1}_{|y_s|\neq 0}\langle
y_s,g(s,y_s,z_s)\rangle\ {\rm d}s
-p\int_t^T |y_s|^{p-2}\mathbbm{1}_{|y_s|\neq 0}\langle y_s,z_s{\rm d}B_s\rangle.
\end{array}
$$
Assumption (A2) yields that, with
probability one, for each $t\in \T$,
$$\begin{array}{lll}
&&\Dis |y_t|^p+c(p)\int_t^T |y_s|^{p-2}\mathbbm{1}_{|y_s|\neq 0}|z_s|^2\ {\rm d}s\\
 &\leq & \Dis|\xi|^p +p\int_t^T
[\psi(|y_s|^p)+\lam |y_s|^{p-1}|z_s|+|y_s|^{p-1}f_s]\ {\rm d}s-p\int_t^T |y_s|^{p-2}\mathbbm{1}_{|y_s|\neq 0}\langle y_s,z_s{\rm d}B_s\rangle.
\end{array}
$$
First of all, in view of the fact that $\psi(\cdot)$ increases at most linearly since it is a nondecreasing concave function and $\psi(0)=0$, we deduce from the previous inequality that
$$\int_0^T |y_s|^{p-2}\mathbbm{1}_{|y_s|\neq 0}|z_s|^2\ {\rm d}s<+\infty,\ \ \ps.$$
Moreover, it follows from the inequality $ab\leq (a^2+b^2)/2 $ that
$$
\begin{array}{lll}
\Dis p\lambda |y_s|^{p-1}|z_s|&=&\Dis p\left({\sqrt{2}\lambda
\over \sqrt{(p-1)\wedge 1}}|y_s|^{p\over 2}\right)
\left(\sqrt{{(p-1)\wedge 1}\over 2}|y_s|^{p-2\over 2}\mathbbm{1}_{|y_s|\neq 0}|z_s|\right)\\
&\leq & \Dis {p\lambda^2\over (p-1)\wedge 1}|y_s|^p +{c(p)\over
2}|y_s|^{p-2}\mathbbm{1}_{|y_s|\neq 0}|z_s|^2.
\end{array}$$
Thus, for each $t\in \T$, we have
\begin{equation}
\Dis |y_t|^p+{c(p)\over 2}\int_t^T |y_s|^{p-2}\mathbbm{1}_{|y_s|\neq
0}|z_s|^2\ {\rm d}s \leq  \Dis X_t-p\int_t^T
|y_s|^{p-2}\mathbbm{1}_{|y_s|\neq 0}\langle y_s,z_s{\rm d}B_s\rangle,
\end{equation}
where
$$X_t=|\xi|^p +d_{\lambda,p} \int_t^T |y_s|^{p}\ {\rm d}s+p\int_t^T \psi(|y_s|^p)\ {\rm d}s+p\int_t^T |y_s|^{p-1}f_s\ {\rm d}s
$$
with $d_{\lambda,p}={p\lambda^2/[(p-1)\wedge 1
]}>0$.

It follows from the BDG inequality that $\{M_t:=\int_0^t
|y_s|^{p-2}\mathbbm{1}_{|y_s|\neq 0}\langle y_s,z_s{\rm
d}B_s\rangle\}_{t\in\T}$ is a uniformly integrable martingale. In
fact, Young's inequality yields
$$
\begin{array}{lll}
\Dis \E\left[\langle M, M \rangle^{1/2}_T\right]&\leq &
\Dis\E\left[\sup\limits_{s\in [0,T]}|y_s|^{p-1}\cdot\left(
\int_0^T|z_s|^2\ {\rm
d}s\right)^{1/2}\right]\\
&= & \Dis \E\left\{\left(\sup\limits_{s\in
[0,T]}|y_s|^{p}\right)^{p-1\over p}\cdot \left[\left(\int_0^T
|z_s|^2\
{\rm d}s\right)^{p/2}\right]^{1\over p}\right\}\\
&\leq &\Dis {(p-1)\over p}\E\left[\sup\limits_{s\in
[0,T]}|y_s|^p\right]+{1\over p}\E\left[\left( \int_0^T|z_s|^2\ {\rm
d}s\right)^{p/2}\right]\\
&<& +\infty.
\end{array}
$$
Thus, for each $0\leq u\leq t\leq T$, taking the conditional mathematical expectation with respect to $\F_u$ in both sides of the inequality (3) yields both
\begin{equation}
{c(p)\over 2}\E\left[\left.\int_t^T |y_s|^{p-2}\mathbbm{1}_{|y_s|\neq
0}|z_s|^2\ {\rm d}s\right|\F_u\right] \leq \E[\left.X_t\right|\F_u]
\end{equation}
and
\begin{equation}
\E\left[\left.\sup\limits_{s\in
[t,T]}|y_s|^{p}\right|\F_u\right]\leq
\E[\left.X_t\right|\F_u]+k_p\E\left[\left.\left(\langle M, M
\rangle_T-\langle M, M \rangle_t\right) ^{1/2}\right|\F_u\right],
\end{equation}
where we have used the BDG inequality in (5), and $k_p$ is a constant depending only on $p$.

On the other hand, it follows from Young's inequality that for each $0\leq u\leq t\leq T$,
$$
\begin{array}{lll}
&&\Dis k_p\E\left[\left.\left(\langle M, M \rangle_T-\langle M, M
\rangle_t\right) ^{1/2}\right|\F_u\right]\\
&\leq & \Dis k_p\E\left[\left.\sup\limits_{s\in
[t,T]}|y_s|^{p/2}\cdot\left( \int_t^T|y_s|^{p-2}\mathbbm{1}_{|y_s|\neq
0}|z_s|^2\ {\rm
d}s\right)^{1/2}\right|\F_u\right]\\
&\leq &\Dis {1\over 2}\E\left[\left.\sup\limits_{s\in
[t,T]}|y_s|^p\right|\F_u\right]+{k_p^2\over 2}\E\left[\left.
\int_t^T|y_s|^{p-2}\mathbbm{1}_{|y_s|\neq 0}|z_s|^2\ {\rm
d}s\right|\F_u\right].
\end{array}
$$
It then follows from inequalities (4) and (5) that there exists a constant $k'_p>0$ depending only on $p$ such that
$$\E\left[\left.\sup\limits_{s\in
[t,T]}|y_s|^{p}\right|\F_u\right]\leq
k'_p\E[\left.X_t\right|\F_u].$$
For each $0\leq u\leq t\leq T$, applying once again Young's inequality we get, with $k''_p$ is
another constant depending only on $p$,
$$\begin{array}{lll}
&&\Dis pk'_p\E\left[\left.\int_t^T |y_s|^{p-1}f_s\ {\rm d}s\right|\F_u\right]\\
&\leq &\Dis pk'_p\E\left[\left.\sup\limits_{s\in
[t,T]}|y_s|^{p-1}\int_t^T f_s\ {\rm d}s\right|\F_u\right]\\
&=& \Dis \E\left[\left.\left({p\over 2(p-1)} \sup\limits_{s\in
[t,T]}|y_s|^{p}\right)^{p-1\over p}\cdot \left[{pk''_p\over
2}\left(\int_t^T f_s\ {\rm
d}s\right)^{p}\right]^{1\over p}\right|\F_u\right]\\
&\leq & \Dis  {1\over 2}\E\left[\left.\sup\limits_{s\in
[t,T]}|y_s|^{p}\right|\F_u\right]+{k''_p\over
2}\E\left[\left.\left(\int_t^T f_s\ {\rm
d}s\right)^{p}\right|\F_u\right],
\end{array}$$
from which we deduce, combing back to the definition of $X_t$, that
$$\begin{array}{lll}
\Dis \E\left[\left.\sup\limits_{s\in
[t,T]}|y_s|^{p}\right|\F_u\right] & \leq &\Dis
2k'_p\E\left[\left.|\xi|^p +d_{\lambda,p} \int_t^T |y_s|^{p}\ {\rm d}s+ p\int_t^T \psi(|y_s|^p)\ {\rm
d}s\right|\F_u\right]\\
&& \Dis +k''_p\E\left[\left.\left(\int_t^T f_s\ {\rm
d}s\right)^{p}\right|\F_u\right].
\end{array}
$$
Letting
$$h_t=\E\left[\left.\sup_{s\in
[t,T]}|y_s|^{p}\right|\F_u\right]$$
in the previous inequality and using Fubini's Theorem and Jensen's inequality yields, in view of the concavity of $\psi(\cdot)$, that for each $0\leq u\leq t\leq
T$,
$$\begin{array}{lll}
\Dis h_t &\leq &\Dis 2k'_p\E[\left.|\xi|^p\right|\F_u]+2pk'_p\int_t^T
\psi(\E[\left.|y_s|^p\right|\F_u])\ {\rm d}s
 +k''_p\E\left[\left.\left(\int_t^T f_s\
{\rm d}s\right)^{p}\right|\F_u\right]\\
&& \Dis +2k'_pd_{\lambda,p}\int_t^T h_s\ {\rm d}s.
\end{array}
$$
Finally, Gronwall's inequality yields that for each $t\in\T$,
$$\begin{array}{l}
\Dis h_t \leq \Dis
e^{2k'_pd_{\lambda,p}(T-t)}\left\{
2k'_p\E[\left.|\xi|^p\right|\F_u]+2pk'_p\int_t^T
\psi(\E[\left.|y_s|^p\right|\F_u])\ {\rm d}s\right.\\
\Dis\hspace*{3.2cm} \left.+k''_p\E\left[\left.\left(\int_t^T f_s\
{\rm d}s\right)^{p}\right|\F_u\right]\right\},
\end{array}
$$
which completes the proof of Proposition 3.\vspace{0.2cm}\hfill $\square$

\section{A stability theorem and the proof of Theorem 1}

In this section, we shall put forward and prove a stability theorem for $L^p\ (p>1)$ solutions to multidimensional BSDEs with generators satisfying (H1)$_{p\wedge 2}$ and (H4). Based on this result, we shall further give the proof of Theorem 1 in Section 2.

The following lemma will be used, which comes from \citet{Fan10}.

{\bf Lemma 1}\ Assume that $\kappa(\cdot):\R^+\mapsto \R^+$ is a nondecreasing and concave function with $\kappa(0)=0$. Then, it increases at most linearly, i.e., there exists a constant $A>0$ such that $$\kappa(x)\leq A(x+1),\ \ \RE\ x\geq 0.$$
Furthermore, for each $m\geq 1$, we have
$$
\kappa(x)\leq (m+2A)x+\kappa({2A\over m+2A}),\ \ \RE\ x\in \R^+.
$$

In the sequel, let $p>1$ and for each $n\geq 1$, let $(y_t,z_t)_{t\in\T}$ and $(y_t^n,z_t^n)_{t\in\T}$ be respectively an $L^p$ solution of the BSDE $(\xi,T,g)$ and the following BSDE depending on parameter $n$:
$$
y_t^n=\xi^n+\int_t^Tg^n (s,y_s^n,z_s^n)
{\rm d}s-\int_t^Tz_s^n {\rm d}B_s,\ \ t\in\T.\vspace{0.2cm}
$$
Furthermore, we introduce the following assumptions:\vspace{0.2cm}

{\bf (B1)}\ \ $\xi^n\in \Lp$ for each $n\geq 1$ and all of $g^n$ satisfy assumptions (H1)$_{p\wedge 2}$ and (H4) with the same $\rho(\cdot)$ and $\bar \lam$.

{\bf (B2)}\ $\Dis \lim\limits_{n\To\infty}\E\left[|\xi^n-\xi|^p+\left(\int_0^T
|g^n(s,y_s,z_s)-g(s,y_s,z_s)|\ {\rm d}s\right)^p\right]=0.$\vspace{0.3cm}

The following Theorem 2 is one of the main results of this section.\vspace{0.2cm}

{\bf Theorem 2}\ \ Under assumptions (B1) and (B2), we have
\begin{equation}
\lim\limits_{n\To\infty}\E\left[\sup\limits_{s\in
[0,T]}|y_s^n-y_s|^p+\left(\int_0^T |z_s^n-z_s|^2\
{\rm d}s\right)^{p/2}\right]=0.
\end{equation}

{\bf Proof.} First, in view of (B1), by $(i)$ of Proposition 1 we note that for each $n\geq 1$,
\begin{itemize}
\item[{\bf (a)}] (H1)$_{p}$ holds true for each $g^n$, together with a new and same function $\hat \rho(x)$;\\ (in case of $1<p\leq 2$, \ \ $\hat \rho(x)\equiv \rho(x)$)
\item[{\bf (b)}] (H1)$_{2}$ holds also true for each $g^n$, together with a new and same function $\bar \rho(x)$.\\ (in case of $p\geq 2$, \ \ $\bar \rho(x)\equiv \rho(x)$)
\end{itemize}

In the sequel, for each $n\geq 1$, let $\hat{y}_\cdot^{n}=y_\cdot^{n}-y_\cdot$,
$\hat{z}_\cdot^{n}=z_\cdot^{n}-z_\cdot$, and
$\hat{\xi}^{n}=\xi^{n}-\xi$. Then
\begin{equation}
\hat{y}_t^{n}=\hat\xi^n+\int_t^T
\hat{g}^{n}(s,\hat{y}_s^{n}, \hat{z}_s^{n})\ {\rm
d}s-\int_t^T \hat{z}_s^{n}{\rm d}B_s,\ \ \ t\in \T,
\end{equation}
where for each $(y,z)\in \R^k\times \R^{k\times d}$,
$$\hat{g}^{n}(s,y,z):=g^n (s,y+y_s,z+z_s)-
g(s,y_s,z_s).$$
Note that
\begin{equation}
\hat{g}^{n}(s,y,z)=g^n (s,y+y_s,z+z_s)-g^n (s,y_s,z_s)
+g^n (s,y_s,z_s)-g(s,y_s,z_s).
\end{equation}
We can check by assumptions (B1) and (B2) together with (a) that the generator ${g}^n$ of BSDE (7) satisfies assumption (A2) with
$$\psi(x)=\hat \rho(x), \ \ \lam=\bar\lam,\ \ {\rm and}\ \ f_t=|g^n (t,y_t,z_t)-g(t,y_t,z_t)|.$$
It then from Proposition 3 with $u=0$ that there exists a constant $C_{\bar\lam,p,T}>0$ depending only on $\bar\lam$, $p$ and $T$ such that for each $n\geq 1$ and each $t\in \T$,
\begin{equation}
\begin{array}{lll}
\Dis\E\left[\sup\limits_{r\in [t,T]}|\hat
y^n_r|^p\right]&\leq & \Dis
C_{\bar\lam,p,T}\E\left[|\hat\xi^n|^p\right]+C_{\bar\lam,p,T}\int_t^T\hat\rho\left(\E\left[\sup_{r\in [s,T]}|\hat y^n_r|^p \right]\right){\rm d}s\\
&& \Dis +C_{\bar\lam,p,T}\E\left[\left(\int_0^T
|g^n (s,y_s,z_s)-g(s,y_s,z_s)|\ {\rm
d}s\right)^p\right].
\end{array}
\end{equation}
Furthermore, in view of (B2) and the fact that $\hat\rho(\cdot)$ is of linear growth by Lemma 1, Gronwall's inequality yields the existence of a constant $M>0$ independent of $n$ such that
$$
\E\left[\sup\limits_{r\in [0,T]}|\hat y^n_r|^p\right]\leq M.
$$
Thus, in view of (B2), by taking the limsup in (9) with respect to $n$ and using Fatou's lemma, the monotonicity and continuity of $\hat\rho(\cdot)$ and Bahari's inequality we can conclude that for each $t\in \T$,
\begin{equation}
\lim\limits_{n\To \infty}\E\left[\sup\limits_{s\in
[t,T]}|y_s^n-y_s|^p\right]= 0.
\end{equation}

Furthermore, by (B2), (8), (b) and Lemma 1 we can also check that the generator ${g}^n$ of BSDE (7) satisfies assumption (A1) with
$$\mu=m+2A,\ \ \lam=\bar\lam,\ \ f_t=|g^n (t,y_t,z_t)-g(t,y_t,z_t)|\ \ {\rm and}\ \ \varphi_t=\bar\rho({2A\over m+2A})$$
for each $m\geq 1$. It then from Proposition 2 with $u=t=0$ that there exists a constant $C_{m,\bar\lam,p,T}>0$ depending on $m$, $\bar\lam$, $p$ and $T$, and a constant $C_{p}$ depending only on $p$ such that for each $m,n\geq 1$,
$$
\begin{array}{lll}
\Dis\E\left[\left(\int_0^T |\hat z^n_s|^2\ {\rm
d}s\right)^{p/2}\right]&\leq & \Dis C_{m,\bar\lam,p,T}\E\left[\sup\limits_{t\in [0,T]}|\hat y^n_t|^p\right]+
C_{p}\left(\bar\rho({2A\over m+2A})\cdot T\right)^{p/2} \\
&& \Dis +C_{p}\E\left[\left(\int_0^T
|g^n (s,y_s,z_s)-g(s,y_s,z_s)|\ {\rm
d}s\right)^p\right].
\end{array}
$$
Thus, in view of (10), (B2) and the fact that $\bar \rho(x)$ is continuous function with $\bar\rho(0)=0$, letting first $n\To \infty$ and then $m\To\infty$ in the previous inequality yields that
$$
\lim\limits_{n\To \infty}\E\left[\left(\int_0^T |z_s^n-z_s|^2\
{\rm d}s\right)^{p/2}\right]= 0.
$$
Thus, we obtain (6). The proof of Theorem 2 is then complete.\vspace{0.2cm}\hfill $\square$

Now, we are in a position to prove Theorem 1.\vspace{0.1cm}

{\bf Proof of Theorem 1.}\ \ Assume that $p>1$, and assumptions (H1)$_{p\wedge 2}$ with $\rho(x)$, (H2)-(H4) and (H5)$_p$ hold for the generator $g$. By $(i)$ of Proposition 1 we note that (H1)$_{2}$ holds also true for $g$, together with a new function $\bar \rho(x)$ (in case of $p\geq 2$, $\bar \rho(x)\equiv \rho(x)$).

The uniqueness part of Theorem 1 is an immediate corollary of Theorem 2. Now, let us prove the existence part. First, for each $n\geq 1$, let $q_n(x)=xn/(|x|\vee n)$ for $x\in\R^k$, and
\begin{equation}
\xi_n:=q_n(\xi)\ \ {\rm and}\ \
g_n(t,y,z):=g(t,y,z)-g(t,0,0)+q_n(g(t,0,0)).
\end{equation}
Note that for each $n\geq 1$, assumptions (H1)$_{2}$ with $\bar \rho(x)$, (H2)-(H4) hold true for each generator $g_n$. Furthermore, for each $n\geq 1$,
\begin{equation}
|\xi_n|\leq n\ \ \ps\ \ {\rm and}\ \  |g_n(t,0,0)|\leq n\ \ \as,
\end{equation}
and by (H5)$_p$ we have
\begin{equation}
\lim\limits_{m,n\To\infty}\E\left[|\xi_{m}-\xi_n|^p+\left(\int_0^T |q_{m}(g(s,0,0))-q_n(g(s,0,0))|\ {\rm d}s\right)^{p}\right]=0.
\end{equation}
By virtue of Theorem 1 in \citet{Xu14} we can know that the BSDE $(\xi_n,T,g_n)$ has a unique $L^2$ solution for each $n\geq 1$, denoted by $(y^n_t,z^n_t)_{t\in \T}$.

Since for each $n\geq 1$, $g_n$ satisfies (H1)$_2$ with $\bar \rho(x)$, and (H4), we can check that it also satisfies (A2) with $$p=2,\ \  \psi(x)=\bar \rho(x), \ \ \lam=\bar\lam \ \ {\rm and}\ \  f_t=q_n(g(t,0,0)).$$
Thus, Proposition 3 together with (12) yields that for each $n\geq 1$, $(y_t^n)_{t\in\T}$ is a bounded process and then belongs to ${\s}^p(0,T;\R^k)$. Furthermore, by Lemma 1 we know that there exists a constant $A>0$ such that
$$\bar \rho(x)\leq A(x+1),\ \ \RE\ x\geq 0, $$
and then $g_n$ satisfies (A1) with
$$\mu=A,\ \ \lam=\bar\lam,\ \ f_t=q_n(g(t,0,0))\ \ {\rm and}\ \  \varphi_t=A,$$
and Proposition 2 together with (12) yields that for each $n\geq 1$, $(z_t^n)_{t\in\T}$ belongs to ${\rm M}^p(0,T;\R^k)$.

In the sequel, for each $m,n\geq 1$, let $$\hat\xi^{m,n}=\xi_{m}-\xi_n,\ \  \hat{y}_\cdot^{m,n}=y_\cdot^{m}-y_\cdot^{n},\ \
\hat{z}_\cdot^{m,n}=z_\cdot^{m}-z_\cdot^{n}.$$
Then $(\hat{y}_\cdot^{m,n},\hat{z}_\cdot^{m,n})$ is an $L^p$ solution of the following BSDE depending on $(m,n)$:
\begin{equation}
\hat{y}_t^{m,n}=\hat\xi^{m,n}+\int_t^T
\hat{g}^{m,n}(s,\hat{y}_s^{m,n}, \hat{z}_s^{m,n})\ {\rm d}s-\int_t^T \hat{z}_s^{m,n}{\rm d}B_s,\ \ \ t\in \T,
\end{equation}
where for each $(y,z)\in \R^k\times \R^{k\times d}$,
$$\hat{g}^{m,n}(s,y,z):=g_{m}(s,y+y_s^n,z+z_s^n)-g_n(s,y_s^n,z_s^n).$$
Note by (11) that for each $m,n\geq 1$,
$$
\hat{g}^{m,n}(t,y,z)=q_{m}(g(t,0,0))-q_n(g(t,0,0))+g(t,y+y_t^n,z+z_t^n)
-g(t,y_t^n,z_t^n).
$$
By the assumptions of the generator $g$ together with (13) we can check that the generator $\hat {g}^{m,n}$ of BSDE (14) satisfies (H1)$_{p\wedge 2}$ and (H4) with $\rho(\cdot)$ and $\bar\lam$ for each $m,n\geq 1$, and
$$
\lim\limits_{m,n\To\infty}\E\left[|\hat \xi^{m,n}-0|^p+\left(\int_0^T
|\hat g^{m,n}(s,0,0)-\tilde{g}(s,0,0)|\ {\rm d}s\right)^p\right]=0,
$$
where for each $(y,z)\in\R^k\times\R^{k\times d}$,
$$\tilde{g}(s,y,z):=0.$$
Thus, we can apply Theorem 2 for BSDE (14) to get that
$$
\lim\limits_{m,n\To\infty}\E\left[\sup\limits_{s\in [0,T]}|\hat y^{m,n}_s-0|^p+\left(\int_0^T |\hat z^{m,n}_s-0|^2\ {\rm
d}s\right)^{p/2}\right]=0,
$$
which means that $\{(y^n_t,z^n_t)_{t\in \T}\}_{n=1}^{\infty}$ is a Cauchy sequence in ${\s}^p(0,T;\R^k)\times {\rm M}^p(0,T;\R^{k\times d})$.

Finally, let $(y_t,z_t)_{t\in \T}$ be the limit process of the sequence $\{(y^n_t,z^n_t)_{t\in \T}\}_{n=1}^{\infty}$ in the process space ${\s}^p(0,T;\R^k)\times{\rm M}^p(0,T;\R^{k\times d})$. We pass to the limit in the sense of uniform convergence in probability for BSDE $(\xi_n,T,g_n)$, thanks to (H2), (H3) and (H4), to see that $(y_t,z_t)_{t\in \T}$ solves the BSDE $(\xi,T,g)$. Thus, we prove the existence part and finally complete the proof of Theorem 1.\hfill $\square$

\section{A comparison theorem}

In this section, we restrict ourselves to the case $k=1$ and prove the following comparison theorem of $L^p$ solutions for BSDEs with generators satisfying (H1)$_p$ and (H4).

{\bf Theorem 3} Let $p>1$, $\xi,\xi'\in \Lp$, $g$ and $g'$ be two generators of BSDEs, and $ (y_\cdot,z_\cdot)$ and $(y'_\cdot,z'_\cdot)$) be respectively an $L^p$ solution to the BSDE $(\xi,T,g)$ and BSDE $(\xi',T,g')$. If $\xi\leq \xi',\ \ps$ and one of the following two statements holds true:
\begin{itemize}
\item[(i)] $g$ satisfies (H1)$_p$ and (H4), and
$$g(t,y'_t,z'_t)\leq g'(t,y'_t,z'_t),\ \ \as;$$
\item[(ii)] $g'$ satisfies (H1)$_p$ and (H4), and $$g(t,y_t,z_t)\leq g'(t,y_t,z_t),\ \ \as,$$
\end{itemize}
then for each $t\in\T$, we have\vspace{0.1cm}
$$y_t\leq y'_t,\ \ \ps.$$

{\bf Proof.}\ \ We first assume that $\xi\leq \xi',\ \ps$, $g$ satisfies (H1)$_p$ with $\rho(x)$ and (H4), and $g(t,y'_t,z'_t)\leq g'(t,y'_t,z'_t),\ \as$.\vspace{0.1cm}

Setting $\hat{y}_t=y_t-y'_t,\ \hat{z}_t=z_t-z'_t,\ \hat \xi=\xi-\xi'$, we have that for each $t\in\T$,
\begin{equation}
\begin{array}{ll}
&\Dis (\hat y_t^+)^p+c(p)\int_t^T |\hat y_s|^{p-2}\mathbbm{1}_{\hat y_s>0}|\hat z_s|^2\ {\rm d}s \\
\leq & \Dis (\hat \xi^+)^p+p\int_t^T |\hat y_s|^{p-1}\mathbbm{1}_{\hat y_s> 0}[g(s,y_s,z_s)-g'(s,y'_s,z'_s)]\ {\rm d}s-p\int_t^T |\hat y_s|^{p-1}\mathbbm{1}_{\hat y_s> 0}\hat z_s{\rm d}B_s
\end{array}
\end{equation}
with $c(p)=p[(p-1)\wedge 1]/2$. Since $g(s,y'_s,z'_s)-g'(s,y'_s,z'_s)$ is non-positive, we
have
\begin{align*}
g(s,y_s,z_s)-g'(s,y'_s,z'_s)&=g(s,y_s,z_s)-g(s,y'_s,z'_s)
+g(s,y'_s,z'_s)-g'(s,y'_s,z'_s)\\
&\leq g(s,y_s,z_s)-g(s,y'_s,z_s)+g(s,y'_s,z_s)-g(s,y'_s,z'_s)
\end{align*}
and we deduce, using (H1)$_p$ and (H4) for $g$ together with a similar inequality before (3), that
\begin{equation}
\begin{array}{lll}
p|\hat y_s|^{p-1}\mathbbm{1}_{\hat y_s> 0}[g(s,y_s,z_s)-g'(s,y'_s,z'_s)]&\leq & \Dis p\rho((\hat{y}_s^+)^p)+p\bar\lam|\hat{y}_s^+|^{p-1}|\hat{z}_s|\vspace{0.1cm}\\
&\leq & \Dis p\bar\rho((\hat{y}_s^+)^p)+{c(p)\over 2}|\hat y_s|^{p-2}\mathbbm{1}_{\hat y_s>0}|\hat z_s|^2,
\end{array}
\end{equation}
where
$$\bar\rho(u):=\rho(u)+d_{\bar\lam,p}u$$
with $d_{\bar\lam,p}=\bar \lam^2/[(p-1)\wedge 1]$ is again a nondecreasing and concave function with $\bar\rho(0)=0$ and $\bar\rho(u)>0$ for $u>0$. Thus, in view of $\xi\leq \xi'$, it follows from (15) and (16) that for each $t\in\T$,
$$
(\hat{y}_t^+)^p \leq p\int_t^T \bar\rho((\hat{y}_s^+)^p){\rm
d}s-p\int_t^T |\hat y_s|^{p-1}\mathbbm{1}_{\hat y_s> 0}\hat z_s{\rm d}B_s,$$
Note that $(\int_0^t |\hat y_s|^{p-1}\mathbbm{1}_{\hat y_s> 0}\hat z_s{\rm d} B_s)_{t\in\T}$ is a martingale by the BDG inequality, and $\bar\rho(u)$ is a concave function. It follows from Jensen's inequality that for each $t\in \T$,
\begin{equation}
\E[(\hat{y}_t^+)^p] \leq p\int_t^T
\bar\rho(\E[(\hat{y}_s^+)^p]){\rm d}s.
\end{equation}
Furthermore, since $\rho(\cdot)$ is a concave function and
$\rho(0)=0$, we have $\rho(u)\geq \rho(1)u$ for each $u\in [0,1]$, and then for each $0\leq u\leq 1$,
$${1\over \bar\rho(u)}={1\over \rho(u)+d_{\bar\lam,p} u}\geq
{1\over \rho(u)+{d_{\bar\lam,p}\over \rho(1)}\rho(u)}={\rho(1)\over
\rho(1)+d_{\bar\lam,p}}\cdot {1\over \rho(u)}.$$
As a result,
$$\int_{0^+}{{\rm d}u\over \bar\rho(u)}=+\infty.$$
Thus, in view of (17), Bihari's inequality yields that for each $t\in \T$, $$\E[(\hat{y}_t^+)^p]=0$$
and then
$$y_t\leq y'_t,\ \ \ps.$$

Now, let us assume that $\xi\leq \xi',\ \ps$, $g'$ satisfies (H1)$_p$ with $\rho(x)$, and (H4), and $ g(t,y_t,z_t)\leq g'(t,y_t,z_t),\ \as$. Since
$g(s,y_s,z_s)-g'(s,y_s,z_s)$ is non-positive, we have
\begin{align*}
g(s,y_s,z_s)-g'(s,y'_s,z'_s)&= g(s,y_s,z_s)-g'(s,y_s,z_s)+g'(s,y_s,z_s)-g'(s,y'_s,z'_s)\\
&\leq g'(s,y_s,z_s)-g'(s,y'_s,z_s)+g'(s,y'_s,z_s)-g'(s,y'_s,z'_s).
\end{align*}
Furthermore, using (H1)$_p$ and (H4) for $g'$, we know that the inequality (16) holds still true. Therefore, the same argument as
above yields that for each $t\in\T$,
$$y_t\leq y'_t,\ \ \ps.$$
The theorem is proved.\vspace{0.2cm}\hfill $\square$

From Theorem 3, the following corollary is immediate.

{\bf Corollary 4} \ Assume that $p>1$ and one of $g$ and $g'$ satisfies (H1)$_p$ and (H4). Let $\xi,\xi'\in \Lp$, and let $(y_\cdot,z_\cdot)$ and $(y'_\cdot,z'_\cdot)$ be respectively an $L^p$ solution to the BSDE $(\xi,T,g)$ and BSDE $(\xi',T,g')$. If $\xi\leq \xi',\ \ps$ and
$$\RE\ y,z,\ \ g(t,y,z)\leq g'(t,y,z)\ \ \as,$$
then for each $t\in\T$,
$$y_t\leq y'_t,\ \ \ps.$$

\section*{Appendix}

In this section, we will give the proof of Proposition 1 in Section 2. The following Lemma 2 will be used frequently, which comes from Lemma 6.1 in \citet{Fan13}.

{\bf Lemma 2}\ Let $f(\cdot)$ be a nondecreasing continuous function on $\R^+$ with $f(0)=0$. Then, the following two statements hold true:

(a) If $f(\cdot)$ is concave on $\R^+$, then $f(x)/x,\  x>0$ is a non-increasing function.

(b) If $f(x)/x,\ x>0$ is a non-increasing function on $\R^+$, then there exists a nondecreasing concave function $p(\cdot)$ defined on $\R^+$ such that for each $x\geq 0$,
$$f(x)\leq p(x) \leq 2f(x).$$

{\bf Proof of Proposition 1}\  First, Let us prove that for each $1\leq p<q$,
$${\rm (H1)}_p\Longrightarrow{\rm (H1)}_q.$$
In fact, assume that $g$ satisfies (H1)$_p$ with a nondecreasing concave function $\rho(\cdot)$. Then we have, $\as$, $\RE y_1,y_2\in\R^k,z\in\R^{k\times d}$,
$$|y_1-y_2|^{q-1}\langle {y_1-y_2\over |y_1-y_2|}\mathbbm{1}_{|y_1-y_2|\neq 0},g(\omega,t,y_1,z)-g(\omega,t,y_2,z)\rangle\leq \bar \rho(|y_1-y_2|^q).$$
where for each $x\geq 0$,
$$\bar \rho(x)=x^{1-{p\over q}}\rho(x^{p\over q}).$$ Obviously, $\bar \rho(\cdot)$ is a nondecreasing continuous function with $\bar\rho(0)=0$ and $\bar \rho(x)>0$ for $x>0$. It follows from (a) of Lemma 2 that
$${\bar\rho(x)\over x}={\rho(x^{p\over q}) \over x^{p\over q}}$$
is a non-increasing function on $\R^+$. Furthermore, by virtue of (b) of Lemma 2 we know that there exists a nondecreasing concave function $\kappa(\cdot)$ such that for each $x\geq 0$,
$$\bar\rho(x)\leq \kappa(x)\leq 2\bar\rho(x).$$
Then we have, $\as$, $\RE y_1,y_2\in\R^k,z\in\R^{k\times d}$,
$$|y_1-y_2|^{q-1}\langle {y_1-y_2\over |y_1-y_2|}\mathbbm{1}_{|y_1-y_2|\neq 0},g(\omega,t,y_1,z)-g(\omega,t,y_2,z)\rangle\leq \kappa(|y_1-y_2|^q)$$
and
$$\int_{0^+}{{\rm d}u\over \kappa(u)}\geq {1\over 2}\int_{0^+} {{\rm d}u\over \bar\rho(u)}={1\over 2}\int_{0^+} {u^{{p\over q}-1}\over \rho(u^{p\over q})}{\rm d}u={q\over 2p}\int_{0^+}{{\rm d}x\over\rho(x)}=+\infty.\vspace{0.2cm}$$
Hence, $g$ satisfies (H1)$_q$ with $\kappa(\cdot)$, and then $(i)$ of Proposition 1 holds true.\vspace{0.2cm}

Then, we prove that for each $1\leq p<q$,
$${\rm (H1b)}_q\Longrightarrow{\rm (H1b)}_p.$$
Indeed, it suffice to show that for a nondecreasing concave function $\rho(\cdot)$ on $\R^+$ with $\rho(0)=0$, if
$$\int_{0+}{u^{q-1}\over \rho^q(u)}{\rm d}u=+\infty,$$
then
$$\int_{0+}{u^{p-1}\over \rho^p(u)}{\rm d}u=+\infty.\vspace{0.1cm}
$$
However, by (a) of Lemma 2 this statement follows easily from the following observation:
$$
\liminf_{u\To 0^+}{{ u^{p-1}\over \rho^p(u)}
\over { u^{q-1}\over \rho^q(u)}}=\liminf_{u\To 0^+}
\left({\rho(u)\over u}\right)^{q-p}\geq \left({\rho(1)\over 1}\right)^{q-p}>0.\vspace{0.1cm}
$$
Hence, $(ii)$ of Proposition 1 is also true.\vspace{0.2cm}

In the sequel, we prove that for each $p\geq 1$,
$${\rm (H1a)}_p \Longrightarrow {\rm (H1b)}_p.$$
In fact, assume that $g$ satisfies (H1a)$_p$ with a nondecreasing concave function $\rho(\cdot)$. Then we have, $\as$, $\RE y_1,y_2\in\R^k,z\in\R^{k\times d}$,
$$\langle {y_1-y_2\over |y_1-y_2|}\mathbbm{1}_{|y_1-y_2|\neq 0},g(\omega,t,y_1,z)-g(\omega,t,y_2,z)\rangle\leq \bar \rho(|y_1-y_2|).$$
where for each $x\geq 0$,
$$\bar \rho(x)=\rho^{1\over p}(x^p).$$
Obviously, $\bar \rho(\cdot)$ is a nondecreasing continuous function with $\bar\rho(0)=0$ and $\bar \rho(x)>0$ for $x>0$. It follows from (a) of Lemma 2 that
$${\bar\rho(x)\over x}=\left({\rho(x^p) \over x^p}\right)^{1\over p}$$
is a non-increasing function on $\R^+$. Furthermore, by virtue of (b) of Lemma 2 we know that there exists a nondecreasing concave function $\kappa(\cdot)$ such that for each $x\geq 0$,
$$\bar\rho(x)\leq \kappa(x)\leq 2\bar\rho(x).$$
Then we have, $\as$, $\RE y_1,y_2\in\R^k,z\in\R^{k\times d}$,
$$\langle {y_1-y_2\over |y_1-y_2|}\mathbbm{1}_{|y_1-y_2|\neq 0},g(\omega,t,y_1,z)-g(\omega,t,y_2,z)\rangle\leq \kappa(|y_1-y_2|)$$
and
$$\int_{0^+}{u^{p-1}\over \kappa^p(u)}{\rm d}u\geq {1\over 2^p}\int_{0^+}{u^{p-1}\over \bar\rho^p(u)}{\rm d}u={1\over 2^p}\int_{0^+} {u^{p-1}\over \rho(u^p)}{\rm d}u={1\over p2^p}\int_{0^+}{1\over\rho(x)}{\rm d}x=+\infty.\vspace{0.1cm}$$
Hence, $g$ satisfies (H1b)$_p$ with $\kappa(\cdot)$.\vspace{0.2cm}

Furthermore, we prove that for each $p\geq 1$,
$${\rm (H1b)}_p \Longrightarrow {\rm (H1a)}_p.$$
In fact, assume that $g$ satisfies (H1b)$_p$ with a nondecreasing concave function $\rho(\cdot)$. Then we have, $\as$, $\RE y_1,y_2\in\R^k,z\in\R^{k\times d}$,
$$\langle {y_1-y_2\over |y_1-y_2|}\mathbbm{1}_{|y_1-y_2|\neq 0},g(\omega,t,y_1,z)-g(\omega,t,y_2,z)\rangle\leq \bar \rho^{1\over p}(|y_1-y_2|^p).$$
where for each $x\geq 0$,
$$\bar \rho(x)=\rho^p(x^{1\over p}).$$
Obviously, $\bar \rho(\cdot)$ is a nondecreasing continuous function with $\bar\rho(0)=0$ and $\bar \rho(x)>0$ for $x>0$. It follows from (a) of Lemma 2 that
$${\bar\rho(x)\over x}=\left({\rho(x^{1\over p}) \over x^{1\over p}}\right)^p$$
is a non-increasing function on $\R^+$. Furthermore, by virtue of (b) of Lemma 2 we know that there exists a nondecreasing concave function $\kappa(\cdot)$ such that for each $x\geq 0$,
$$\bar\rho(x)\leq \kappa(x)\leq 2\bar\rho(x).$$
Then we have, $\as$, $\RE y_1,y_2\in\R^k,z\in\R^{k\times d}$,
$$\langle {y_1-y_2\over |y_1-y_2|}\mathbbm{1}_{|y_1-y_2|\neq 0},g(\omega,t,y_1,z)-g(\omega,t,y_2,z)\rangle\leq \kappa^{1\over p}(|y_1-y_2|^p)$$
and
$$\int_{0^+}{{\rm d}u\over \kappa(u)}\geq {1\over 2}\int_{0^+} {{\rm d}u\over \bar\rho(u)}={1\over 2}\int_{0^+} {{\rm d}u\over \rho^p(u^{1\over p})}={p\over 2}\int_{0^+}{x^{p-1}\over\rho^p(x)}{\rm d}x=+\infty.\vspace{0.1cm}$$
Hence, $g$ satisfies (H1a)$_p$ with $\kappa(\cdot)$, and then $(iii)$ of Proposition 1 holds true.\vspace{0.2cm}

Finally, we prove that the concavity condition of $\rho(\cdot)$ in (H1b)$_p$ and (H1a)$_p$ can be replaced with the continuity condition.

Assume first that $p\geq 1$ and $g$ satisfies (H1b)$_p$ with $\rho(\cdot):\R^+\mapsto\R^+$, which is a nondecreasing continuous (but not concave) function with $\rho(0)=0$, $\rho(u)>0$ for $u>0$ and
$$\int_{0^+} {u^{p-1}\over \rho^p(u)}{\rm d}u=+\infty.\vspace{0.1cm}$$
Then, there exists a ${\mathcal P}\subset \Omega\times \T$ with
$${\rm d}{\mathbbm P}\times {\rm d}t((\Omega\times \T)\cap {\mathcal P}^c)=0$$
such that for each $(\omega,t)\in {\mathcal P}$, $y_1,y_2\in \R^k$ and $z\in\R^{k\times d}$, we have
\begin{equation}
\langle {y_1-y_2\over |y_1-y_2|}\mathbbm{1}_{|y_1-y_2|\neq 0},g(\omega,t,y_1,z)-g(\omega,t,y_2,z)\rangle\leq \rho(|y_1-y_2|).
\end{equation}
Now, for each $r\in \R^+$, let
$$
\begin{array}{l}
F(r)=\sup\left\{\langle {y_1-y_2\over |y_1-y_2|}\mathbbm{1}_{|y_1-y_2|\neq 0},g(\omega,t,y_1,z)-g(\omega,t,y_2,z)\rangle
:\ (y_1,y_2)\in\R^k\times\R^k,\right. \\
\hspace{3.5cm}\Dis |y_1-y_2|\leq r,\ \left.(\omega,t,z)\in {\mathcal P}\times \R^{k\times d}\right\}.
\end{array}
$$
It is clear that $F(0)=0$. It follows from (18) that $F(\cdot)$ is well-defined, nondecreasing and for each $r\geq 0$,
$$0\leq F(r)\leq \rho(r).$$
Thus, in view of the continuity of $\rho(\cdot)$ and the fact $\rho(0)=0$, we know that $F(\cdot)$ is right-continuous at $0$. Furthermore, for $r,s\geq 0$ and $(y_1,y_2)\in \R^k\times \R^k$ with
$r\leq |y_1-y_2|\leq r+s$, it follows from the definition of $F(\cdot)$ that for each $(\omega,t)\in {\mathcal P}$, $y_1,y_2\in \R^k$ and $z\in\R^{k\times d}$,
$$
\begin{array}{lll}
&&\Dis\langle {y_1-y_2\over |y_1-y_2|}\mathbbm{1}_{|y_1-y_2|\neq 0},g(\omega,t,y_1,z)-g(\omega,t,y_2,z)\rangle\\
&=&\Dis \langle {y_1-y_2\over |y_1-y_2|}\mathbbm{1}_{|y_1-y_2|\neq 0},g(\omega,t,y_1,z)-g(\omega,t,y_1+r{y_2-y_1\over |y_2-y_1|},z)\rangle\\
&&\Dis +\langle {y_1-y_2\over |y_1-y_2|}\mathbbm{1}_{|y_1-y_2|\neq 0},g(\omega,t,y_1+r{y_2-y_1\over |y_2-y_1|},z)-g(\omega,t,y_2,z)\rangle\vspace{0.2cm}\\
&\leq&\Dis F(r)+F(s)
\end{array}
$$
so that, the case $|y_1-y_2|\leq r$ being trivial, we conclude that $F(\cdot)$ is sub-additive. That is, for each $r,s\geq 0$,
$$F(r+s)\leq F(r)+F(s).$$
As a result, $F(\cdot)$ is a continuous modular function, and then there exists a nondecreasing and concave function $\bar\rho(\cdot):\R^+\To \R^+$ such that for each $u\geq 0$,
$$F(u)\leq \bar\rho(u)\leq 2F(u)\leq 2\rho(u)$$
(see pages 499-500 in \cite{Kuang04} for details). Thus, it follows from (18) and the definition of $F(\cdot)$ that $\as$, for each $(y_1,y_2,z)\in \R^k\times \R^k\times \R^{k\times d}$,
$$\langle {y_1-y_2\over |y_1-y_2|}\mathbbm{1}_{|y_1-y_2|\neq 0},g(\omega,t,y_1,z)-g(\omega,t,y_2,z)\rangle\leq F(|y_1-y_2|)\leq \bar\rho(|y_1-y_2|)\leq \kappa(|y_1-y_2|),$$
where
$$\kappa(u):=\bar\rho(u)+u.$$
Clearly, $\kappa(\cdot)$ is a nondecreasing and concave function with $\kappa(0)=0$ and $\kappa(u)>0$ for $u>0$. Thus, for completing the proof that $g$ satisfies (H1b)$_p$, it suffices to show that
$$\int_{0^+} {u^{p-1}\over \kappa^p(u)}{\rm d}u=+\infty.$$
Indeed, if $\bar\rho(u)>0$ for each $u>0$, since $\bar\rho(\cdot)$ is concave on $\R^+$ with $\bar\rho(0)=0$, we have $\bar\rho(u)\geq u\bar\rho(1)$ for each $u\in (0,1)$, and then
$$
\begin{array}{lll}
\Dis\int_{0^+}{u^{p-1}\over \kappa^p(u)}{\rm d}u&=&\Dis\int_{0^+} {u^{p-1}\over (\bar\rho(u)+u)^p}{\rm d}u\vspace{0.1cm}\\
&\geq &\Dis {\bar\rho^p(1)\over (1+\bar\rho(1))^p}\int_{0^+}{u^{p-1}\over \bar\rho^p(u)}{\rm d}u\vspace{0.1cm}\\
&\geq &\Dis {1\over 2^p}\cdot {\bar\rho^p(1)\over (1+\bar\rho(1))^p}\int_{0^+} {u^{p-1}\over \rho^p(u)}{\rm d}u\vspace{0.1cm}\\
&=& \Dis +\infty.
\end{array}
$$
Otherwise,
$$\int_{0^+} {u^{p-1}\over \kappa^p(u)}{\rm d}u=\int_{0^+} {u^{p-1}\over (\bar\rho(u)+u)^p}{\rm d}u=\int_{0^+} {{\rm d}u\over u}=+\infty.\vspace{0.2cm}$$
Thus, in view of $(iii)$ of Proposition 1, we have proved that the concavity condition of $\rho(\cdot)$ in (H1b)$_p$ and (H1a)$_p$ can be replaced with the continuity condition. The proof of Proposition 1 is then completed.\hfill $\square$\vspace{0.2cm}




\begin{thebibliography}{00}

\setlength{\baselineskip}{13.5pt}
\setlength{\itemsep}{1.8mm}

\bibitem[Bahlali(2002)]{Bah02}Bahlali, K., 2002. Existence and uniqueness of solutions for BSDEs with locally Lipschitz coefficient. Electronic Communications in probability 7:169--179.

\bibitem[Bahlali, Essaky and Hassani(2010)]{Bah10}Bahlali, K., Essaky, E., Hassani, M., 2010. Multidimensional BSDEs with super-linear growth coefficient: Application to degenerate systems of semilinear PDEs. C. R. Acad. Sci. Paris, Ser. I 348:677--682.

\bibitem[Bihari(1956)]{Bih56}Bihari, I., 1956. A generalization of a lemma of Bellman and its application to uniqueness problem of differential equations. Acta Math. Acad. Sci. Hungar. 7: 71--94.

\bibitem[Briand and Confortola(2008)]{BC08}Briand, Ph., Confortola, F., 2008. BSDEs with stochastic Lipschitz condition and quadratic PDEs in Hilbert spaces. Stochastic Processes and Their Applications  118:818--838.

\bibitem[Briand et al.(2003)]{Bri03}Briand, Ph., Delyon, B., Hu, Y., Pardoux, E., Stoica, L., 2003. $L^p$ solutions of backward stochastic differential equations. Stochastic Processes and Their Applications 108:109--129.

\bibitem[Briand and Hu(2006)]{Bri06}Briand, Ph., Hu, Y., 2006. BSDE with quadratic growth and unbounded terminal value. Probability Theory and Related Fields 136:604--618.

\bibitem[Briand and Hu(2008)]{Bri08}Briand, Ph., Hu, Y., 2008. Quadratic BSDEs with convex generators and unbounded terminal conditions. Probability Theory and Related Fields 141:543--567.

\bibitem[Briand, Lepeltier and San Martin(2007)]{Bri07}Briand, Ph., Lepeltier, J.-P., San Martin, J., 2007. One-dimensional BSDEs whose coefficient is monotonic in $y$ and non-Lipschitz in $z$. Bernoulli 13:80--91.

\bibitem[Buckdahn, Hu and Peng(1999)]{Buc99}Buckdahn, R., Hu, Y., Peng, S., 1999. Probabilistic approach to homogenization of viscosity solutions of parabolic PDEs. NoDEA Nonlinear Differential Equations Appl. 6(4):395--411.

\bibitem[Buckdahn, Quincampoix and Rascanu(2000)]{Buc00}Buckdahn, R., Quincampoix, M., Rascanu, A., 2000. Viability property for a backward stochastic differential equations and applications to partial differential equations. Probability Theory and Related Fields 116:485--504.

\bibitem[Chen(2010)]{Chen10}Chen, S., 2010. $L^p$ solutions of one-dimensional backward stochastic differential equations with continuous coefficients. Stochastic Analysis and Applications 28:820--841.

\bibitem[Chen and Wang(2000)]{Chen00}Chen, Z., Wang, B., 2000. Infinite time interval BSDEs and the convergence of $g$-martingales. J. Austral. Math. Soc. (Series A) 69:187--211.

\bibitem[Constantin(2001)]{Cons01}Constantin, G., 2001. On the existence and uniqueness of adapted solutions for backward stochastic differential equations. Analele Universit\u{a}\c{t}ii din Timi\c{s}oara,Seria Matematic\u{a}-Informatic\u{a} XXXIX(2):15--22.

\bibitem[Delbaen, Hu and Bao(2011)]{Delb11} Delbaen, F., Hu, Y.,  Bao, X., 2011. Backward SDEs with superquadratic growth. Probability Theory and Related Fields 150(24):145--192.

\bibitem[Delbaen and Tang(2010)]{Delb10} Delbaen, F., Tang, S., 2010. Harmonic analysis of stochastic equations and backward stochastic differential equations. Probability Theory and Related Fields 146:291--336.

\bibitem[El Karoui, Peng and Quenez(1997)]{El97}El Karoui, N., Peng, S., Quenez, M.C., 1997. Backward stochastic differential equations in finance. Math. Finance 7:1--72.

\bibitem[Fan and Jiang(2010)]{Fan10}Fan, S., Jiang, L., 2010. Uniqueness result for the BSDE whose generator is monotonic in y and uniformly continuous in z. C. R. Acad. Sci. Paris, Ser. I 348:89--92.

\bibitem[Fan and Jiang(2011)]{FJ11}Fan, S., Jiang, L., 2011. Existence and uniqueness result for a backward stochastic differential equation whose generator is Lipschitz continuous in $y$ and uniformly continuous in $z$. Journal of Applied Mathematics and Computing 36:1--10.

\bibitem[Fan and Jiang(2012)]{Fan12} Fan, S., Jiang, L., 2012. One-dimensional BSDEs with left-continuous, lower semi-continuous and linear-growth generators. Statistics and Probability Letters 82:1792--1798.

\bibitem[Fan and Jiang(2012)]{FJ12a}Fan, S., Jiang, L., 2012. A generalized comparison theorem for BSDEs and its applications. Journal of Theoretical Probability 25(1):50--61.

\bibitem[Fan and Jiang(2012)]{FJ12b}Fan, S., Jiang, L., 2012. $L^p\ (p>1)$ solutions for one-dimensional BSDEs with linear-growth generators. Journal of Applied Mathematics and Computing 38(1-2):295--304.

\bibitem[Fan and Jiang(2013)]{Fan13}Fan, S., Jiang, L., 2013. Multidimensional BSDEs with weakly monotonic generators. Acta Mathematica Sinica, English Series 23(10):1885--1906.

\bibitem[Fan and Jiang(2014)]{Fan14}Fan, S., Jiang, L., 2014. $L^p$ solutions of BSDEs with a new kind of non-Lipschitz coefficients. To appear in Acta Mathematics Applicatae Sinica, English Series. arXiv: 1402.6773v1 [math.PR].

\bibitem[Fan, Jiang and Davison(2010)]{Fan10b} Fan, S., Jiang, L., Davison, M., 2010. Uniqueness of solutions for multidimensional BSDEs with uniformly continuous generators, C. R. Acad. Sci. Paris, Ser. I, 348: 683--686.

\bibitem[Fan, Jiang and Tian(2011)]{Fan11} Fan, S., Jiang, L., Tian, D., 2011. One-dimensional BSDEs with finite and infinite time horizons. Stochastic Processes and Their Applications 121:427--440.

\bibitem[Hamad\`{e}ne(2003)]{Ham03}Hamad\`{e}ne, S., 2003. Multidimensional backward stochastic differential equations with uniformly continuous coefficients. Bernoulli 9(3):517--534.

\bibitem[Jia(2008)]{Jia08}Jia, G., 2008. A class of backward stochastic differential equations with discontinuous coefficients. Statistics and Probability Letters 78:231--237.

\bibitem[Jia(2010)]{Jia10}Jia, G., 2010. Backward Stochastic differential equations with a uniformly continuous generator and related g-expectation. Stochastic Processes and Their Applications 120(11):2241--2257.

\bibitem[Kobylanski(2000)]{Kob00} Kobylanski, M., 2000. Backward stochastic differential equations and partial equations with quadratic growth. Annals of Probability 28:259--276.

\bibitem[Kuang(2004)]{Kuang04}Kuang, J., 2004. Applied Inequalities, 3nd edition, Shandong Science and Technology Press, Jinan (in Chinese).

\bibitem[Lepeltier and San Martin(1997)]{Lep97}Lepeltier, J.-P., San Martin, J., 1997. Backward stochastic differential equations with continuous coefficient. Statistics and Probability Letters 32:425--430.

\bibitem[Lepeltier and San Mart\'{\i}n(1998)]{Lep98} Lepeltier, J.-P., San Mart\'{\i}n, J., 1998. Existence for BSDE with superlinear quadratic coefficient. Stochastics and Stochastic Reports 63:227--240.

\bibitem[Mao(1995)]{Mao95}Mao, X., 1995. Adapted solutions of backward stochastic differential
    equations with non-Lipschitz coefficients. Stochastic Process and Their Applications 58:281-292.

\bibitem[Morlais(2009)]{Mor09}Morlais, M.-A, 2009. Quadratic BSDEs driven by a continuous martingale and applications to the utility maximization problem. Finance Stoch. 13:121--150.

\bibitem[Pardoux and Peng(1990)]{Par90}Pardoux, E., Peng, S., 1990. Adapted solution of a backward stochastic differential equation. Systems Control Letters 14:55--61.

\bibitem[Pardoux(1999)]{Par99}Pardoux, E., 1999. BSDEs, weak convergence and homogenization of semilinear PDEs. Nonlinear Analysis, Differential Equations and Control (Montreal, QC, 1998). Kluwer Academic Publishers, Dordrecht, pp.503--549.

\bibitem[Peng(1991)]{Peng91}Peng, S., 1991. Probabilistic interpretation for systems of quasilinear parabolic partial differential equations. Stochatics Stochatics Reports 37:61--74.

\bibitem[Peng(1997)]{Peng97}Peng, S., 1997. Backward SDE and related g-expectation. In: El Karoui, N., Mazliak, L. (Eds.), Backward Stochastic Differential Equations, Pitman Research Notes Mathematical Series, Vol. 364 Longman, Harlow, pp.141--159.

\bibitem[Situ(1997)]{Situ97}Situ, R., 1997. On the solutions of backward stochastic differential equations with jumps and applications. Stochastic Process and Their Applications 66:209--236.

\bibitem[Tang and Li(1994)]{Tang94} Tang, S., Li, X., 1994. Necessary conditions for optimal control of stochastic systems with random jumps. SIAM J. Control Optim. 32(5):1447--1475

\bibitem[Tang(1998)]{Tang98}Tang, S., 1998. The maximum principle for partially observed optimal control of stochastic differential equations[J]. SIAM J. Control Optim. 36(5):1596-1617.

\bibitem[Xu and Fan(2014)]{Xu14}Xu, S., Fan, S., 2014. A general existence and uniqueness result on mutidimensional BSDEs. arXiv: 1402.6777v1 [math.PR].

\bibitem[Yao(2010)]{Yao10}Yao, S., 2010. $L^p$ solutions for backward stochastic differential equations with jumps. arXiv: 1007.2226v1 [math.PR].

\bibitem[Yin and Mao(2008)]{Yin08}Yin, J., Mao, X., 2008. The adapted solution and comparison theorem for backward stochastic differential equations with poisson jumps and applications. Journal of Mathematical Analysis and Applications 346:345--358.

\bibitem[Yin and Situ(2003)]{Yin03}Yin, J., Situ, R., 2003. On solutions of foward-backward stochastic differential equation with Possion jumps. Stochastic Analysis and Applications 23(6):1419--1448.

\end{thebibliography}
\end{document}